\documentclass{svmult}
\usepackage{amsmath}
\usepackage[english]{babel}
\usepackage{amstext}
\usepackage{amsmath}
\usepackage{amsfonts}
\usepackage[all]{xy}
\UseTips
\usepackage{url}

\spnewtheorem{thm}{Theorem}{\bfseries}{\itshape}
\spnewtheorem{defn}[thm]{Definition}{\bfseries}{\rmfamily}
\spnewtheorem{rmk}[thm]{Remark}{\bfseries}{\rmfamily}
\spnewtheorem{rmks}[thm]{Remarks}{\bfseries}{\rmfamily}
\spnewtheorem{prop}[thm]{Proposition}{\bfseries}{\itshape}
\spnewtheorem{ex}[thm]{Example}{\bfseries}{\rmfamily}
\numberwithin{thm}{section}

\newcommand{\p}{\partial}

\newcommand{\eq}[2]{\begin{equation}\label{#1}#2 \end{equation}}

\newcommand{\ga}[2]{\begin{gather}\label{#1}#2 \end{gather}}

\newcommand{\Pic}{{\rm Pic}}
\newcommand{\Hom}{{\rm Hom}}
\newcommand{\Spec}{{\rm Spec \,}}

\newcommand{\Char}{{\rm char\,}}
\newcommand{\Tr}{{\rm Tr}}

\newcommand{\Div}{{\rm div}}
\newcommand{\bw}[2]{\mathbb{W}_{#1}({#2})}
\newcommand{\bwc}[3]{\mathbb{W}_{#1}\Omega^{#2}_{#3}}
\newcommand{\w}[2]{\textrm{W}_{#1}({#2})}
\newcommand{\wc}[3]{\textrm{W}_{#1}\Omega^{#2}_{#3}}
\newcommand{\pp}{\mathbb{P}^1}
\newcommand{\Th}[3]{\mbox{$\textrm{TH}^{#1}(#2,{#1};{#3})$}}
\newcommand{\Ch}{\textrm{CH}}
\newcommand{\Sh}{\textrm{SH}}

\newcommand{\wV}{\textrm{V}}
\newcommand{\wF}{\textrm{F}}
\newcommand{\wR}{\textrm{R}}

\newcommand{\Zy}{\textrm{Z}}

\newcommand{\sL}{{\mathcal L}}
\newcommand{\sM}{{\mathcal M}}

\newcommand{\sO}{{\mathcal O}}

\newcommand{\A}{{\mathbb A}}

\newcommand{\C}{{\mathbb C}}
\newcommand{\F}{{\mathbb F}}

\newcommand{\N}{{\mathbb N}}
\renewcommand{\P}{{\mathbb P}}
\newcommand{\Q}{{\mathbb Q}}
\newcommand{\R}{{\mathbb R}}

\newcommand{\Z}{{\mathbb Z}}

\begin{document}

\title*{Characteristic 0 and $p$ analogies, and some motivic cohomology}
\titlerunning{Characteristic 0 and $p$, and motivic cohomology}
\author{Manuel Blickle, H\'el\`ene Esnault, Kay R\"ulling\thanks{This work has been partly supported by the DFG Schwerpunkt ``Komplexe Mannigfaltigkeiten'' und bei the DFG Leibniz Program}}
\institute{Universit\"at Duisburg-Essen, Mathematik, 45117 Essen, Germany
\texttt{manuel.blickle@uni-essen.de}\\
\texttt{esnault@uni-essen.de}\\
\texttt{kay.ruelling@uni-essen.de}}

\maketitle
\section*{Introduction}
The purpose of this survey is to explain some recent results about analogies between
characteristic 0 and characteristic $p>0$ geometry, and to discuss an infinitesimal
variant of motivic cohomology.

Homotopy invariance for motivic cohomology implies, in particular, that the Picard
group of the affine line over a filed $k$ is trivial, i.e. ${\rm Pic}(\A^1_k)=0$.
However, if instead of considering the Picard group, we consider the group of
isomorphism classes of pairs $(\sL, t)$ consisting of a line bundle $\sL$ on
$\A^1_k$,
 and an isomorphism $t: \sO_{{\rm Spec}(\sO_{\A^1}/\frak{m}^n)}\xrightarrow{\,{\cong}\,} \sL_{{\rm
Spec}(\sO_{\A^1}/\frak{m}^n)}$ where $\frak{m}$ is the maximal ideal at the origin,
then one obtains the group of 1 units in $k[[t]]/(t^n)$, that is $\{\text{isom.
classes }(\sL, t) \}\cong \big(1+tk[[t]]/(1+t^nk[[t]])\big)^\times$. This group is
indeed a ring, namely the big ring of Witt vectors $\bw{n-1}{k}$ of length $n-1$
over $k$.

On the other hand, the theory of additive higher Chow groups, or Chow groups with
modulus condition 2, developed in \cite{BlEs03a} and \cite{BlEs03b}, allows to
realize absolute differential forms $\Omega^i_{k/\Z}$  of $k$  as a group of
0-cycles (Theorems \ref{2.1} and \ref{2.3}). In section 2 these groups of 0-cycles
 and their higher modulus $n$ generalization are presented. Theorem \ref{2.5} --
the main result of section 2 -- asserts that these groups of 0-cycles with higher
modulus $n$ compute big Witt differential forms $\bwc{n-1}{i}{k}$. Cutting out the
$p$-isotypical component in characteristic $p>0$ by suitable correspondences, and
extending Theorem \ref{2.5} to  smooth local rings over a perfect field would then
describe crystalline cohomology of a smooth proper variety over a perfect field of
characteristic $p>0$  as the hypercohomology of a complex of sheaves of 0-cycles.

If $X$ is a smooth complex variety, the Riemann-Hilbert correspondence establishes
an equivalence of categories between holonomic $D_X$-modules (coherent) and
constructible sheaves (topological). Emerton and Kisin developed in
\cite{EmKis.Fcrys} for characteristic $p>0$ a correspondence with properties
analogous to the complex theory. The ``topological like'' side consists of sheaves
of finite dimensional $\F_p$-vector spaces which are locally constant in the \'etale
topology, the ``$D$-module like'' side consists of locally finitely generated unit
$\sO_{X,F}$ modules (see Theorem \ref{thm.main}), which are modules on which the
Frobenius acts in a certain way (unit). Even if the unit condition is quite
restrictive -- as will be explaind shortly, it essentially corresponds to ``slope
zero'' -- the theory has a vaste range of applications and many objects naturally
carry such a structure. For example $\sO_X$ itself is such a unit $\sO_{X,F}$
module. In section 3 we investigate singularities of $Y$ (at a point $x \in Y
\subseteq X$, where $X$ is smooth) from this viewpoint and obtain striking analogies
for local invariants over the complex numbers and in characteristic $p>0$.
Postponing the definition of these invariants, which were introduced by Lyubeznik in
\cite{Lyub.FinChar0}, the main result (Theorem \ref{thm.main}) of section 3 asserts
that these invariants can be expressed in terms of \'etale cohomology
$H^i_{x}(Y_{\operatorname{\acute{e}t}}, \F_p)$ in characteristic $p>0$ and
respectively in terms of singular cohomology $H^i_{x}(Y_{\operatorname{an}},
\mathbb{C})$ over the complex numbers, by virtually the same expression. This
extends earlier results in the isolated singular analytic case of
\cite{Lyub.FinChar0} and \cite{LopezSabbah}. This extension is made possible by the
similarity between the two correspondences which allows to essentially treat both
settings (char. $p>0$ and char. 0) formally as one.

In section 1, we review recent results on congruences modulo $q$-powers for the
number of $\F_q$-rational points of algebraic varieties. They are all based on
Deligne's philosophy which predicts a deep analogy between the level of congruences
for varieties defined over $\F_q$ and the Hodge level for varieties  over the
complex numbers. If a variety has Hodge level $\ge \kappa$ over the complex numbers,
one expects that ``over'' $\F_q$ it will have the same number of rational points as
$\P^n$ modulo $q^\kappa$. Of course, the challenge is to make precise ``over'' as it
can't be the same variety. Divisibility of eigenvalues of the geometric Frobenius
acting on $\ell$-adic cohomology is one method to show the existence of congruences,
slope computation in crystalline cohomology, or, in the singular case, in
Berthelot's rigid cohomology is another one. Of course, ideally one would wish to
prove a motivic statement which would imply all those results at once. But it is
often beyond reach. The main results are Theorem \ref{thm1.5}, which in particular
gives a positive answer to the Lang-Manin conjecture asserting that Fano varieties
over a finite field have a rational point, Theorem \ref{thm1.7}, asserting that the
mod $p$ reduction of a regular model of a smooth projective variety defined over a
local field, the $\ell$-adic cohomology of which is supported in codimension 1,
carries one rational point modulo $q$, Theorem \ref{thm1.11}, asserting that two
theta divisors on an abelian variety defined over a finite field carry the same
number of points modulo $q$. Theorem \ref{thm1.11} answers positively the finite
consequence of a conjecture of Serre, and appears as a consequence of the slope
Theorem \ref{thm1.10} asserting that the slope $<1$ piece of rigid cohomology is
computed by Witt vector cohomology. Theorem \ref{thm1.5} can be proven using either
$\ell$-adic cohomology or crystalline cohomology. Indeed, the geometric result
behind is that Fano varieties are rationally connected and therefore their Chow
group of 0-cycles satisfies base change. Theorem \ref{thm1.7} relies on a
generalization to local fields of Deligne's integrality theorem over finite fields
(Theorem \ref{thm1.8}).\\ \ \\
\section{Hodge type over the complex numbers and congruences for the number of rational points over a finite field}
\subsection{Deligne's integrality theorem over a finite field}\label{subsec1.1}
Deligne developed the  theory of weights for complex varieties via the weight
filtration in his mixed Hodge theory \cite{DHodge}, and, for varieties defined over
finite fields, via the  absolute values with respect to any complex embedding of the
eigenvalues of the geometric Frobenius operator acting on $\ell$-adic cohomology
\cite{DWeil}. The philosophy of motives, as conceived by him and Grothendieck,
predicts a closed analogy between those two concepts of weights. It has led to many
very fundamental results, the first of which being Deligne's proof of the Weil
conjecture for Hodge level one  complete intersections and for $K3$ surfaces
(Invent. math. {\bf 15} (1972)).
 On the other hand, Deligne shows the fundamental integrality theorem\\ \ \\
\begin{thm}[Deligne \cite{Dint}, Corollaire 5.5.3]\label{thm1.1}  Let $X$ be a scheme of finite type defined over $\F_q$. Then  the eigenvalues of the geometric Frobenius acting on  compactly supported $\ell$-adic cohomology $H^i_c(X\times_{\F_q} \overline{\F_q}, \Q_\ell)$ are algebraic integers.
\end{thm}
If $X$ is defined over the complex numbers, its Hodge filtration $F^j$ satisfies
${\rm gr}_F^jH^i_c(X)=0$ for $j<0$ (see  Hodge III \cite{DHodge}). In Deligne's  motivic
philosophy, the integrality theorem \ref{thm1.1} is analogous to the Hodge
filtration starting in degrees $\ge 0$ in  Hodge theory. This analogy has been less
studied in the past than the weight analogy. The purpose of this section is to
demonstrate on some examples how the analogy works  between the $F$-filtration in
Hodge theory and the integrality over a finite field for $\ell$-adic cohomology. It
has led in very recent years to
a series of results on congruences for the number of points on varieties defined over finite fields. \\ \ \\
Theorem \ref{thm1.1} implies
\begin{thm} \label{thm1.2} Let $X$ be a scheme of finite type defined over $\F_q$. Then  the eigenvalues of the geometric Frobenius acting on   $\ell$-adic cohomology $H^i(X\times_{\F_q} \overline{\F_q}, \Q_\ell)$ are algebraic integers.
\end{thm}
Strictly speaking, Deligne shows this for $X$ smooth via duality, but applying de
Jong's alteration theorem, one easily reduces the theorem to the smooth case as in
\cite{DE}, Corollary 0.3.
\subsection{Divisibility and rational points} \label{subsec1.2}
The Grothendieck-Lefschetz trace formula \cite{Gr} \ga{1}{
|X(\F_q)|=\sum_{i=0}^\infty (-1)^i {\rm Trace \ Frobenius}|H^i_c(X\times_{\F_q}
\overline{\F_q}, \Q_\ell) } together with Theorem \ref{thm1.1} implies that if the
eigenvalues of the geometric Frobenius are not only algebraic integers, i.e. $\in
\bar{\Z}$, but also, for $i\ge 1$,  they are $q$-divisible as algebraic integers,
i.e. $\in q\cdot\bar{\Z}$,  then  one has \ga{2}{|X(\F_q)|\equiv {\rm dim}\
H^0_c(X\times_{\F_q}\overline{\F_q}) \ {\rm mod}\ q.} So we conclude \ga{3}{
 X \ {\rm  proper \ and \ geometrically \  connected}\\
{\rm with \ eigenvalues \ of \ geom. \ Frob.   \ } \in q\cdot \bar{\Z} \notag\\
 \Longrightarrow
|X(\F_q)|\equiv 1 \ {\rm mod} \ q. \notag} Purity,  for which  the smoothness
condition is definitely necessary, together with Theorem \ref{thm1.2}  implies
\begin{thm} \label{thm1.3} Let $X$ be a smooth scheme of finite type defined over $\F_q$. Then  the eigenvalues of the geometric Frobenius acting on   $\ell$-adic cohomology $H^i_{A\times_{\F_q} \overline{\F_q}}(X\times_{\F_q} \overline{\F_q}, \Q_\ell)$ with supports along $A$ are algebraic integers divisible by $q^\kappa$, if $A\subset X$ is a closed subscheme of codimension $\ge \kappa$.
\end{thm}
(See \cite{EFano}, Lemma 2.1 for the analogous proof in crystalline cohomology). If
$X$ is no longer smooth, the conclusion of Theorem \ref{thm1.3} is no longer true
(see \cite{DE}, Remark 0.5). But it remains true for generic hyperplanes for
example, as purity is then true (see \cite{Eeigv}, Theorem 2.1). In Deligne's
philosophy, Theorem \ref{thm1.3} is analogous to the Hodge level statement
\begin{thm} \label{thm1.4}
Let $X$ be a smooth scheme of finite type defined over $\C$. Then the graded pieces
for the Hodge filtration  $F$ on de Rham cohomology with support fulfill
 $ {\rm gr}_F^a H^i_{A}(X)=0$ for $a <\kappa$, if $A\subset X$ is a closed subscheme of codimension $\ge \kappa$.
\end{thm}
The proof, based on purity, is the same as for Theorem \ref{thm1.3}.
\subsection{Chow group of 0-cycles, coniveau and divisibility for smooth proper varieties} \label{subsec1.3}
So if $X$ is smooth proper geometrically connected over $\F_q$, Theorem \ref{thm1.3}
implies that \eqref{3}  holds if $H^i(X\times_{\F_q} \overline{\F_q}, \Q_\ell)$ is
supported in codimension $\ge 1$ for $i\ge 1$. An equivalent terminology is to say
that  $H^i(X\times_{\F_q} \overline{\F_q}, \Q_\ell)$ has coniveau $\ge 1$. According
to Bloch's decomposition of the diagonal (\cite{B}, Appendix to Lecture 1), this is
the case if the Chow group of $0$-cycles over a field containing the field of
rational functions $\F_q(X)$ is trivial. We conclude
\begin{thm}[\cite{EFano}, Corollary 1.2, Corollary 1.3] \label{thm1.5} Let $X$ be a smooth, proper,  geometrically connected variety over $\F_q$. Assume that  $CH_0(X\times_{\F_q} \overline{\F_q(X)})=\Q$.  Then $|X(\F_q)|\equiv 1 $ mod $q$. In particular, Fano varieties have a rational point, as conjectured by Lang \cite{L} and Manin \cite{Ma}.
\end{thm}
Originally, Bloch decomposed the diagonal in order to recover Mumford's theorem (and
its variants) asserting in its simplest form that if $X$ is a smooth projective
surface over $\C$, and if ${\rm gr}_F^0H^2(X)\neq 0$, then it can't be true that
$CH_0(X)_{{\rm deg}=0}\cong {\rm Alb}(X)$. Furthermore, that Fano varieties, that is
smooth projective geometrically irreducible varieties $X$ so that the inverse of the
dualizing sheaf $\omega_X^{-1}$ is ample, are rationally connected over any
algebraically closed field, is a consequence of Mori's break and bend theory, and
has been proven independently by Koll\'ar-Miyaoka-Mori and Campana (see \cite{Ko}
and references there).
\subsection{Singular varieties defined by equations: Ax-Katz' theorem and divisibility} \label{subsec1.4}
If $X$ is no longer smooth, not only one can't apply Theorem \ref{thm1.3} to get
divisibility of eigenvalues, but also Bloch's decomposition of the diagonal does not
work. Indeed, the diagonal has only a homology class,  so it does not act on
cohomology, while Grothendieck-Lefschetz trace formula \eqref{1} allowing to count
points needs cohomology.   Yet, Deligne's philosophy on the analogy between Hodge
type over $\C$ and eigenvalue divisibility over $\F_q$, is still at disposal. There
is an instance where one can directly generalize the Leitfaden sketched for the
proof of Theorem \ref{thm1.5}, by refining the motivic cohomology used. Rather than
considering $CH_0(X)$ of the singular variety, one embedds $X\subset \P^n$ in a
projective space, and considers relative motivic cohomology $H^{2n}(\P\times U,
Y\times U, n)$ (\cite{BEL}), section 1), where $\P\xrightarrow{\pi} \P^n$ is an
alteration so that $Y=\pi^{-1}(X)$ is a normal crossings divisor and
$U=\P^n\setminus X$. Then the graph of $\pi$ has a cycle class in this relative
motivic group and one shows (\cite{BEL}, Theorem 1.2) that if $X$ is a hypersurface
in $\P^n$ of degree $\le n$, then this class decomposes in a suitable sense,
generalizing Bloch's decomposability notion in $H^{2{\rm dim}(X)}(X, {\rm
dim}(X))\cong CH_0(X)$ when $X$ is smooth. This implies immediately Hodge type $\ge
1$ over $\C$ as well as
 eigenvalue $q$-divisibility over $\F_q$, and yields a motivic proof of Ax' theorem asserting
the congruence
\eqref{3} for hypersurfaces of degree $\le n$ in $\P^n$ over $\F_q$. Ax' theorem  generalizes the mod $p$ congruence due to Chevalley-Warning. \\ \ \\
Another instance for which one can make Deligne's philosophy work concerns closed
subsets $X$ of $\P^n$ defined by equations of degrees $d_1\ge \ldots \ge d_r$,
without any other assumption. Those equations could be chosen in a highly
non-optimal way. For example, one could take $r$ times the same equation. Ax and
Katz in \cite{Ax}, \cite{Ka}, assign to $X$ a level $\kappa:={\rm max}\{0,
[\frac{n-d_2-\ldots -d_r}{d_1}]\}$, and show that $|X(\F_q)|\equiv |\P^n(\F_q)|$ mod
$q^\kappa$.  On the other hand, one can compute (\cite{EHodgetyp} and \cite{ENS})
that over $\C$, ${\rm gr}_F^aH_c^i(\P^n\setminus X)=0$ for $a<\kappa$. Thus one
expects eigenvalue divisibility. Indeed one has
\begin{thm}[\cite{Eeigv}, Theorem 1.1, \cite{EK}, Theorem 2.1] \label{thm1.6} Let $X\subset \P^n$ be a closed subset defined by equations of degrees $d_1\ge \ldots \ge d_r$. Then the eigenvalues of the geometric Frobenius acting on $H^i_c((\P^n\setminus X)\times_{\F_q}
\overline{\F_q}, \Q_\ell)$ are in $q^\kappa\cdot \bar{\Z}$.
\end{thm}
According to \eqref{1}, Theorem \ref{thm1.6} implies Ax-Katz' theorem. However, the
proof of Theorem \ref{thm1.6} is not good as it uses Ax-Katz' theorem. One would
like to understand a motivic proof in the spirit of \cite{BEL}. We are very far from
it, as, even if $X$ is a smooth hypersurface of low degree,  we do not know how to
compute its Chow groups of higher dimensional cycles.
\subsection{Singular varieties in families} \label{subsec1.5}
Singular varieties which are degenerations of smooth ones have more structure.
Fakhruddin and Rajan (\cite{FR}, Corollary 1.2)   generalize  the motivic method of
Theorem \ref{thm1.5} in a relative situation: if $f: X\to S$ is a proper dominant
morphism of smooth irreducible varieties over a finite field $k$ with
$CH_0(X\times_S \overline{k(X)})=\Q$, then for any $s \in S(k)$, one has
$|f^{-1}(s)| \equiv 1$ mod $|k|$. Similarly on the Hodge side one proves
(\cite{EappFakh}, Theorem 1.1)  that if $f: X\to S$ is a proper morphism with $S$ a
smooth connected curve and $X$ smooth, then if ${\rm gr}_F^0H^i(f^{-1}(s_0))=0$ for
some $s_0$ in the smooth locus of $f$ and all $i\ge 1$, then $ {\rm
gr}_F^0H^i(f^{-1}(s))=0$ for all $s$ and all $i\ge 1$. Those two statements, the
first one in equal characteristic $p>0$ with its  strong motivic assumption, the
second one in equal characteristic 0 with its minimal Hodge type assumption,
suggest, using Deligne's philosophy,  that the mod $p$ reduction of a smooth
projective variety in characteristic zero with $ {\rm gr}_F^0H^i(X)=0$ for all $i\ge
1$ has eigenvalue $q$-divisibility, and therefore by \eqref{1},  its number of
rational points is congruent to one modulo $q$. One shows
\begin{thm}[\cite{Eunqchar}, Theorem 1.1, Section 4, Proof of Theorem 1.1]
 \label{thm1.7} Let $X$ be a smooth projective variety over a local field
$K$ with finite residue field $\F_q$. Assume $H^i(X\times_K \bar{K}, \Q_\ell)$ lives
in codimension $\ge 1$ for all $i\ge 1$. Then the eigenvalues of the geometric
Frobenius acting on $H^i(Y\times_{\F_q} \overline{\F_q}, \Q_\ell) $ for $i\ge 1$,
where $Y$ is the mod $p$ reduction  of a projective regular model, are lying in
$q\cdot \bar{\Z}$. In particular, $|Y(\F_q)|\equiv 1$ mod $q$.
\end{thm}
If the local field $K$  has equal characteristic $p>0$, one  expects Theorem
\ref{thm1.7} to be optimal. However, if the local field $K$ has unequal
characteristic, one would wish, following Deligne's philosophy, to replace the
assumption on the coniveau of $\ell$-adic cohomology  by ${\rm gr}_F^0 H^i(X)=0$ for
all $i\ge 1$. Due to the comparison of \'etale and de Rham cohomology, and due to
the Hodge conjecture in codimension 1, those  two assumptions are equivalent for
surfaces. In general,  Grothendieck's generalized Hodge conjecture in codimension 1
predicts that if $ {\rm gr}_F^0 H^i(X)=0$ then $H^i(X)$ lives in codimension $\ge
1$. Thus those two conditions are expected to be equivalent in general. So in
unequal characteristic, Theorem \ref{thm1.7} is optimal for surfaces, but in higher
dimension,  in absence of a proof of the generalized Hodge conjecture in codimension
1, one would wish to have another proof.

 On the other hand, a generalization of Bloch's decomposition of the diagonal implies that the motivic assumption
$CH_0(X_0\otimes_{K_0} \Omega)=\Q$, where $K_0\subset K$ is a subfield of finite
type over the prime field over which $X$ is defined, i.e. $X=X_0\times_{K_0}K$, and
$\Omega$ contains $K_0(X_0)$, so for example, $\Omega=\bar{K}$ in unequal
characteristic, implies the coniveau assumption of the theorem. So the motivic
assumption implies the following direct corollary of  Theorem \ref{thm1.7}: Let $X$
be a smooth projective variety over a local field $K$ with finite residue field
$\F_q$. Assume $CH_0(X_0\times_{K_0} \Omega)=\Q$. Then the eigenvalues of the
geometric Frobenius acting on $H^i(Y\times_{\F_q} \overline{\F_q}, \Q_\ell) $ for
$i\ge 1$, where $Y$ is the mod $p$ reduction  of a projective regular model, are
lying in $q\cdot \bar{\Z}$. In particular, $|Y(\F_q)|\equiv 1$ mod $q$
(\cite{Eunqchar}, Corollary 1.2).

 It is to be noted that for surfaces in characteristic zero, we are very far from knowing a positive answer to Bloch's conjecture, asserting that  ${\rm gr}_F^0 H^i(X)=0$ for $i\ge 1$ implies
the motivic condition. Thus the range of applicability of Theorem \ref{thm1.7} is
much larger than the one of its corollary.

The proof of Theorem \ref{thm1.7} relies on the specialization map, on Gabber's
purity theorem  \cite{Fu}, Theorem 2.1.1, de Jong's alteration \cite{dJ}, Theorem
6.5, and on the direct generalization of Deligne's integrality theorem \ref{thm1.1}
to the local field case
\begin{thm}[\cite{DE}, Theorem 0.2, Corollary 0.3] \label{thm1.8} Let $X$ be a scheme of finite type defined over a local field $K$ with finite residue field. Then the eigenvalue of a lifting of the geometric Frobenius in the Galois group ${\rm Gal}(\bar{K}/K)$
of the local field acting on $H^i_c(X\times_K \bar{K})$ and $H^i(X\times_K \bar{K})$
are algebraic integers for all $i$.
\end{thm}
For the proof of Theorem \ref{thm1.7},  one needs a   form of the integrality
theorem over local fields which is weaker than the one stated in Theorem
\ref{thm1.8}, and which can be proven directly using alterations. Theorem
\ref{thm1.8} itself is a corollary of a more general integrality theorem for
$\ell$-adic sheaves.

As a corollary, one has, as in Theorem \ref{thm1.3} over finite field, the
divisibility statement for smooth varieties
\begin{thm}[\cite{DE}, Corollary 0.4] \label{thm1.9 }
 Let $X$ be a smooth scheme of finite type defined over a local field $K$ with finite residue field $\F_q$. Then  the eigenvalues of
a lifting of the geometric Frobenius in the Galois group ${\rm Gal}(\bar{K}/K)$ of
the local field acting on $H^i_{A\times_{K} \bar{K}}(X\times_{K} \bar{K}, \Q_\ell)$
are algebraic integers divisible by $q^\kappa$, if $A\subset X$ is a closed
subscheme of codimension $\ge \kappa$.
\end{thm}
\subsection{Witt vector cohomology} \label{subsec.16}
Intuitively, if instead of  $\ell$-adic cohomology, we  consider crystalline or
rigid cohomology, we expect a more direct link between Hodge type and slopes, and
consequently congruences for the number of points if the ground field is finite. The
theorem of Bloch and Illusie asserts that if $X$ is smooth proper over a perfect
field $k$
 of characteristic $p > 0$, $W = W(k)$, $K =
{\rm Frac}(W)$, then  there is a functorial isomorphism (\cite{Bl}, III, 3.5 and
\cite{I}, II, 3.5) \ga{4}{H^{i}(X/K)^{<1}\ \xrightarrow{\cong }\ H^{i}(X,
W\sO_X)_K.} Here $H^{i}(X/K)^{<1}$ denotes the maximal subspace of crystalline
cohomology on which Frobenius acts with slopes $< 1$, the subscript $_K$ denotes
tensorisation with $K$, and the right hand side is Witt vector cohomology, as
considered by Serre in \cite{Se}. On the other hand, over $\C$, one has a functorial
surjective map (\cite{EHodgetyp},  Proposition 1.2, \cite{EappFakh}, Proof of
Theorem 1.1) \ga{5}{ H^{i}(X,\sO_X)\xrightarrow{{\rm surj}}{\rm gr}_0^F H^{i}(X). }
This gives an upper bound on the Hodge type in the sense that if the left hand side,
which is usually easy to compute as this is coherent cohomology, dies, then so does
the right hand side, which is a topological invariant. A weak version of \eqref{5}
allows one for example to remark  that if $\Theta\subset A$ is a theta divisor on an
abelian variety $A$ of dimension $g$ over $\C$ (i.e. effective, ample, with $h^0(A,
\sO_A(\Theta))=1$), then one has an isomorphism \ga{6}{{\rm gr}_0^F H^g_c(A\setminus
\Theta)\xrightarrow{{\rm iso}} {\rm gr}_0^F H^g(A). } There is a generalization of
\eqref{4} and \eqref{5}
\begin{thm}[\cite{BBE}, Theorem 1.1] \label{thm1.10}
Let $X$ be a proper scheme over a perfect field $k$ of characteristic $p>0$, then
\eqref{4} holds true, where the left hand side is replaced by Berthelot's rigid
cohomology. So Witt vector cohomology computes the slope $<1$ piece of rigid
cohomology.
\end{thm}
Then the finite field version of \eqref{6} asserts
\begin{thm}[\cite{BBE}, Theorem 1.4] \label{thm1.11} Let $A$ be an abelian variety defined over $\F_q$, and let $\Theta, \Theta'$ be two theta divisors. Then
$|\Theta(\F_q)|\equiv |\Theta'(\F_q)|$ mod $q$.
\end{thm}
This answers positively the finite field consequence of a  conjecture of Serre
asserting that the difference of the motives of $\Theta$ and $\Theta'$ should be
divisible by the Lefschetz motive over any field.
\\ \ \\
\begin{rmks} \label{rmks1.12}
As concluding remarks to this section, let us first observe that the motivic
philosophy discussed here leads to various questions anchored directly in geometry.
As an example, as already mentioned in \cite{EFano}, section 3, Gorenstein Fano
varieties $X$ in characteristic 0 fulfill ${\rm gr}_F^0 H^i(X)=0$ for $i\ge 1$. This
suggests the existence of a good definition  of rational singularities in
characteristic $p>0$ so that a Gorenstein Fano variety over $\F_q$ would have one
rational point modulo $q$. This would also require a generalization of Mori's
 break and bend method to varieties with this type of mild singularities.

Next, let us remark that we did not discuss in this survey higher congruences, that
is congruences modulo $q^\kappa, \kappa\ge 2$. Indeed, Theorem \ref{thm1.3} implies
that if $X$ is projective smooth over $\F_q$, then if $H^i(X\times_{\F_q}
\overline{\F_q}, \Q_\ell)$ is supported in codimension $\ge \kappa$ up to the class
of the $j$-th self-product  of the polarization if $i=2j$, then \eqref{1} yields
$|X(\F_q)| \equiv |\P^n(\F_q)| \ {\rm mod} \ q^\kappa$. However, following the
Leitfaden explained in subsection \ref{subsec1.3}, the coniveau $\kappa$ condition
is implied by triviality of $CH_i(X\times_{\F_q} \overline{\F_q(X)})$, for all
$i<\kappa$, and this is a condition we can basically never check.

Finally, by Theorem \ref{thm1.10}, Witt vector cohomology  computes the slope $<1$
part of rigid cohomology. In particular, it is a topological cohomology theory. We
do not know the relation between rigid and $\ell$-adic cohomology if $X$ is proper
but not smooth. On the other hand, the slope $=0$ part of crystalline cohomology,
when $X$ is smooth, is easier to understand as it is described by coherent
cohomology.  This viewpoint is developed in section 3.

\end{rmks}

\section{Additive higher Chow groups with higher modulus of type (n,n) over a field}
\subsection{Additive higher Chow groups} \label{subsec2.1}
Let $k$ be a field and $X$ an equidimensional $k$-scheme. In \cite{Bl86} Bloch
develops a theory of higher Chow groups, which are isomorphic to the motivic
cohomology groups $\Ch^p(X,n)\cong H^{2p-n}_\sM(X,\Z(p))$, $p,n\ge 0$ (see
\cite{V02}), by using the scheme  $\Delta^n=\Spec
k[t_0,\ldots,t_n]/\left(\sum_{i=0}^n t_i-1\right)$ or, in the cubical definition,
the scheme $(\pp\setminus\{1\})^n$. Replacing in the simplicial definition $\sum
t_i=1$ by $\sum t_i= \lambda$ yields the same groups as long as $\lambda\in
k^\times$. The degenerate case $\lambda=0$ is investigated in \cite{BlEs03a}. One
obtains a theory of additive higher Chow groups, $\Sh^p(X,n)$, $p\ge 0$, $n\ge 1$.
In analogy to the theorem of Nesterenko-Suslin and Totaro (see \cite{NeSu89},
\cite{To92}) $\Ch^n(k,n)\cong K^M_n(k)$, it is shown
\begin{thm}[\cite{BlEs03a}, Theorem 5.3]\label{2.1}
Let $k$ be field with  $\Char k\not= 2$, then there is an isomorphism of groups
\[\Sh^n(k,n)\cong\Omega^{n-1}_{k/\Z}.\]
\end{thm}
The proof is in the spirit of the proofs in \cite{NeSu89}, \cite{To92}. We will
sketch the proof for the corresponding statement with $\Sh^n(k,n)$ replaced by a
cubical version of the higher additive Chow groups. This cubical version is defined
in \cite{BlEs03b}, so far only for a field and on the level of 0-cycles. The
definition is as follows. Consider the $k$-scheme $X_n=\A^1\times
(\pp\setminus\{1\})^n$ with coordinates $(x,y_1,\ldots,y_n)$ and denote the union of
all faces by $Y_n=\bigcup_{i=1}^n (y_i=0,\infty)$. Now denote by $\Zy_0(k,n-1)$ the
free abelian group generated on all closed points of
$\A^1\setminus\{0\}\times(\pp\setminus\{0,1,\infty\})^{n-1}$. Let $\Zy_1(k,n;2)$ be
the free abelian group generated on all irreducible curves $C\subset X_n\setminus
Y_n$ satisfying the following properties
\begin{enumerate}
\item ({\em Good position}) $\p^j_i[C]=(y_i=j).[C]\in\Zy_0(k,n-1)$, for $i=1,\ldots,n, j=0,\infty$.
\item ({\em Modulus 2 condition}) Let $\nu: \widetilde{C}\to \pp\times(\pp)^n$ be the normalization of the compactification of $C$, then in $\Zy_0(\widetilde{C})$
\ga{7}{
        2\,\Div(\nu^*x)\le \sum_{i=1}^n\Div(\nu^*y_i-1).
      }
\end{enumerate}
Because of (i) one has a complex $$\p=\sum_{i=1}^n(-1)^i (\p^0_i-\p^\infty_i):
\Zy_1(k,n;2)\to\Zy_0(k,n-1)\to 0.$$
\begin{defn}[\cite{BlEs03b}, Definition 6.2]\label{2.2}
The additive higher Chow groups of type $(n,n)$ and modulus 2 of a field $k$ are
given by the homology of the above complex, i.e.
\[\Th{n}{k}{2}=\frac{\Zy_0(k,n-1)}{\p\Zy_1(k,n;2)}.\]
\end{defn}
It is shown in \cite[Theorem 6.4]{BlEs03b} (cf. Theorem \ref{2.1} )
\begin{thm}\label{2.3}
Let $k$ be a field with $\Char k\not= 2,3$, then the map
\eq{7.05}{\Th{n}{k}{2}\stackrel{\simeq}{\longrightarrow}\Omega^{n-1}_{k/\Z},}
\[[P]\mapsto \Tr_{k(P)/k}\left(\frac{1}{x(P)}\frac{dy_1(P)}{y_1(P)}\cdots\frac{dy_{n-1}(P)}{y_{n-1}(P)}\right)\]
is an isomorphism of groups. Furthermore, the inclusion $\iota:
(\pp\setminus\{1\})^{n-1}\hookrightarrow X_{n-1}$, $(y_1,\ldots,y_{n-1})\mapsto
(1,y_1,\ldots,y_{n-1})$ induces a commutative diagram
\ga{}{\xymatrix{K^M_{n-1}(k)\ar[r]^(0.4)\simeq\ar[d]_{\rm dlog} & \Ch^{n-1}(k,n-1)\ar[d]_{\iota_*}\\
            \Omega^{n-1}_{k/\Z}\ar[r]^(0.4)\simeq  & \Th{n}{k}{2}} \notag}
defined on elements by \ga{}{  \xymatrix{
\{a_1,\ldots, a_{n-1}\}\ar@{|->}[r]\ar@{|->}[d] & (a_1,\ldots,a_{n-1})\ar@{|->}[d]\\
\frac{da_1}{a_1}\wedge\ldots\wedge\frac{da_{n-1}}{a_{n-1}}\ar@{|->}[r] & (1,
a_1,\ldots,a_{n-1}).
           }
\notag}
\end{thm}
The idea of the proof is the following. One first shows, that the map \eqref{7.05}
is well defined, using the reciprocity law for rational differential forms on a
non-singular projective curve. Here one explicitly uses the modulus condition
\eqref{7}. Then one constructs an inverse map with the help of a representation of
the absolute differential forms of $k$ by generators and relations. Showing this map
to be well defined is equivalent to finding 1-cycles in $\Zy_1(k,n;2)$, whose
boundary yield the relations. This map is easily seen to be injective and the
surjectivity follows from the fact, that the trace on the absolute differentials
corresponds to the pushforward on the additive Chow groups.

\subsection{Higher modulus} \label{subsec2.2}
The cubical definition allows  one to generalize the definition of the additive
higher Chow groups (of type (n,n)) to higher modulus, i.e. replace the $2$ in
equation \eqref{7} by an integer $m\ge 2$. The resulting groups are denoted by
$\Th{n}{k}{m}$. (Notice that one has
$\Th{n}{k}{0}=\Ch^n(\A^1\setminus\{0\},n-1)=0=\Th{n}{k}{1}$.) Up to now it is not
clear how to formulate this higher modulus condition in the simplicial setup.

The attempt to generalize Theorem \ref{2.3} leads to the following considerations.
Let $\Pic(\A^1, m\{0\})$ be the relative Picard group of $\A^1$ with modulus
$m\{0\}$. Then there is a natural surjective map
\eq{7.1}{\Pic(\A^1,m\{0\})\longrightarrow \Th{1}{k}{m}.} For this one observes that,
if two $0$-cycles $\Div f$ and $\Div g$ in the left hand side are equal, then the
curve $C\subset \A^1\times\pp\setminus\{1\}$ defined by $fy=g$ satisfies the modulus
condition
\[m(x=0).\bar{C}\le (y-1).\bar{C},\]
with $\bar{C}\subset \pp\times\pp$ the closure of $C$. Hence $\Div f=\Div g$ also in
the right hand side. By Theorem \ref{2.3} this is an isomorphism for $m=2$ and one
can show that it is an isomorphism for all $m\ge 2$. On the other hand we may
identify $\Pic(\A^1,m\{0\})$ as a group with the additive group of the ring of big
Witt vectors of length $m-1$ of $k$ via
\eq{7.2}{\Pic(\A^1,m\{0\})\cong\left(\frac{1+tk[t]}{1+t^{m+1}k[t]}\right)^\times\cong\bw{m-1}{k}.}
This together with Theorem \ref{2.3} leads to the prediction, that $\Th{n}{k}{m}$ is
isomorphic to the group of generalized degree $n-1$ Witt differential forms of
length $m-1$. These groups form the generalized de Rham-Witt complex of
Hesselholt-Madsen generalizing the $p$-typical de Rham-Witt complex of
Bloch-Deligne-Illusie (see \cite{Bl}, \cite{I}). Before stating the generalization
of Theorem \ref{2.3} to the case of higher modulus, we describe the de Rham-Witt
complex and some of his properties.
\begin{defn}[see \cite{HeMa01}, \cite{R}, cf. \cite{I}]\label{2.4}
Let $A$ be a ring. A  {\em Witt complex over $A$} is a projective system of
differential graded $\Z$-algebras
\[((E_m)_{m\in\N}, \wR: E_{m+1}\to E_m)\]
together with families of homomorphisms of graded rings
\[(\wF_n: E_{nm+n-1}\to E_{m})_{m,n\in\N}\]
and homomorphisms of graded groups
\[ (V_{n}: E_{m}\to E_{nm+n-1})_{m,n\in\N}\]
satisfying the following relations, for all $n,r\in \N$
\begin{enumerate}
\item[(i)] $\wR\wF_n=\wF_n\wR^n,$ $\wR^n\wV_n=\wV_n\wR$, $\wF_1=\wV_1={\rm id}$, $\wF_n\wF_r=\wF_{nr}$, $\wV_n\wV_r=\wV_{nr}$.
\item[(ii)] $\wF_n\wV_n=n,$ and if $(n,r)=1$, then $\wF_r\wV_n=\wV_n\wF_r$ on $E_{rm+r-1}$.
\item[(iii)] $\wV_n(\wF_n(x)y)=x\wV_n(y)$, for $x\in E_{nm+n-1}$, $y\in E_m$.
\item[(iv)] $\wF_n d \wV_n= d,$ with $d$ the differential on $E_m$ and $E_{nm+n-1}$ respectively.
\end{enumerate}
Furthermore, there is a homomorphism of projective systems of rings
\[(\lambda: \bw{m}{A}\to E_m^0)_{m\in\N},\]
under which the Frobenius and the Verschiebung maps on the big Witt vectors
correspond to $\wF_n$ and $\wV_n$, $n\in\N$ on $E^0$ and satisfies
\begin{enumerate}
\item[(v)] $\wF_n d \lambda([a])= \lambda([a]^{n-1})d\lambda([a]),$ for $a\in A$.
\end{enumerate}
A morphism of Witt complexes over $A$ is a morphism of projective systems of dga's
compatible with all the structures. This yields a category of Witt complexes over
$A$.
\end{defn}
\begin{thm}[see \cite{HeMa01}, \cite{R}, cf. \cite{I} ]
The category of Witt complexes has an initial object, called the de Rham-Witt
complex of $A$ and denoted by $(\bwc{m}{\cdot}{A})_{m\in\\N}$.
\end{thm}
Hesselholt-Madsen prove this using the Freyd adjoint functor theorem. In case $A$ is
a $\Z_{(p)}$-algebra, $p\not=2$ a prime, $\bwc{m}{\cdot}{A}$ may be constructed following
Illusie as a quotient of $\Omega^\cdot_{\bw{m}{A}}$. It follows
\[\bwc{1}{\cdot}{A}=\Omega^\cdot_{A/\Z}\quad \bwc{m}{0}{A}=\bw{m}{A}.\]
\begin{rmk}
Let $A$ be a $\Z_{(p)}$-algebra, $p\not= 2$ a prime, and denote by
$(\wc{n}{\cdot}{A})_{n\in\N_0}$ the $p$-typical de Rham-Witt complex of
Bloch-Deligne-Illusie-Hesselholt-Madsen, then one has
\ga{8}{\bwc{m}{\cdot}{A}\cong\prod_{(j,p)=1}\wc{n(j)}{\cdot}{A}, \quad n(j)\text{ given by }
jp^{n(j)}\le m<j p^{n(j)+1}.} If $X$ is a smooth variety over a perfect field, Bloch
and Illusie also define $\wc{n}{\cdot}{X}$, which is a complex of coherent sheaves on the
scheme $\w{n}{X}$. Its hypercohomology equals the crystalline cohomology. This is
used for example to derive \eqref{4}.
\end{rmk}
The generalization of Theorem \ref{2.3} to the case of higher modulus is given by
the
\begin{thm}[\cite{R}]\label{2.5}
Let $k$ be a field with $\Char\not= 2$. Then the projective system
$(\bigoplus_{n\ge1}\Th{n}{k}{m+1}\to\bigoplus_{n\ge1}\Th{n}{k}{m})_{m\ge 2}$ can be
equipped with the structure of a Witt complex over $k$ and the natural map
\linebreak $\bwc{*}{\cdot}{k}\to\bigoplus_{n\ge 0}\Th{n+1}{k}{*+1}$ induced by the
universality of the de Rham-Witt complex is an isomorphism. In particular
\[\bwc{m-1}{n-1}{k}\cong\Th{n}{k}{m}.\]
\end{thm}
In the following the Witt complex structure of $$T_m:=\linebreak\bigoplus_{n\ge
0}\Th{n+1}{k}{m+1}$$ is described. The multiplication of a graded commutative ring
on $T_m$ is induced by the exterior product of cycles followed by a pushforward,
which is induced by the multiplication map $\A^1\times\A^1\to \A^1$. The
differential is induced by pushing forward via the diagonal $\A^1\to\A^1\times\pp$
and then restricting to $\A^1\times\pp\setminus\{1\}$. $\wF_r$ (resp. $\wV_r$) is
induced by pushing forward (resp. pulling back) via $\A^1\to\A^1$, $a\to a^r$. And
finally the map $\bw{m}{k}\to T_m^0=\Th{1}{k}{m+1}$ is given by the composition of
\eqref{7.2} and \eqref{7.1}. The relations (i)-(iii) and (v) in Definition \ref{2.4}
are already satisfied on the level of cycles. That $T_m$ is a dga with (iv), follows
from the result of Nesterenko-Suslin and Totaro, since one has a surjective map
\[\bigoplus_{P_0\in\A^1\setminus\{0\}}CH^{n-1}(k(P_0),n-1)\to \Th{n}{k}{m}\]
induced by the inclusions
$\{P_0\}\times(\pp\setminus\{1\})^{n-1}\hookrightarrow\A^1\setminus\{0\}\times(\pp\setminus\{1\})^{n-1}$.
Thus $T$ is a Witt complex and this gives the natural map from Theorem \ref{2.5}. In
\cite{R} a trace map for arbitrary field extensions on the de Rham-Witt groups is
constructed as well as a residue symbol for closed points of a smooth projective
curve on the rational Witt differentials of this curve and a reciprocity law is
proven, generalizing the corresponding notions and statements on K\"ahler
differentials (and the one obtained by Witt in \cite{Wi36}). Using this results one
can generalize the proof of Theorem \ref{2.3} to obtain, that the map
\ga{}{\Th{n}{k}{m}\to \bwc{m-1}{n-1}{k} \notag \\
[P]\mapsto
\Tr_{k(P)/k}\left(\frac{1}{[x(P)]}\frac{d[y_1(P)]}{[y_1(P)]}\cdots\frac{d[y_{n-1}(P)]}{[y_{n-1}(P)]}\right)
 \notag}
gives the inverse map to the one obtained by the universality of the de Rham-Witt
complex. Here $[-]: k(P)\to \bw{m-1}{k(P)}$ is the Teichm\"uller lift.

Finally we explain how to describe the $p$-typical de Rham-Witt complex over a field
$k$ via the additive higher Chow groups. Denote by
\[f=\exp(\sum_{i=0}^\infty -\frac{t^{p^i}}{p^i})\in 1+t\Z_{(p)}[[t]]\]
the inverse of the Artin-Hasse eponential and by
\[ f_m\in (1+\oplus_{i=1}^m t^i\Z_{(p)})\]
the truncation of $f$. Write $\epsilon_m$ for the 0-cycle $[\Div
f_m]\in\Th{1}{k}{m+1}$. Then it follows from the description of the additive higher
Chow groups as the de Rham-Witt complex and from \eqref{8}, that one has
\[\wc{r}{n-1}{k}\cong \Th{n}{k}{p^r+1} *\epsilon_{p^r},\]
where we denote by $*$ the multiplication in the additive higher Chow groups
explained above.

%
%

\newcommand{\CC}{\mathbb{C}}
\renewcommand{\AA}{\mathbb{A}}
\newcommand{\FF}{\mathbb{F}}
\newcommand{\QQ}{\mathbb{Q}}
\newcommand{\ZZ}{\mathbb{Z}}
\newcommand{\PP}{\mathbb{P}}
\newcommand{\Oo}{\mathcal{O}}
\newcommand{\Mm}{\mathcal{M}}
\newcommand{\Nn}{\mathcal{N}}
\newcommand{\Yy}{\mathcal{Y}}
\newcommand{\Ll}{\mathcal{L}}
\newcommand{\deR}{\operatorname{dR}}
\newcommand{\tensor}{\otimes}
\renewcommand{\to}[1][]{\xrightarrow{\ #1\ }}
\renewcommand{\Spec}{\operatorname{Spec}}
\newcommand{\op}[1]{\operatorname{#1}}
\newcommand{\Sol}{\operatorname{Sol}}
\newcommand{\defeq}{\stackrel{\scriptscriptstyle \op{def}}{=}}
\renewcommand{\R}{\mathbf{R}}
\newcommand{\RHom}{\op{\mathbf{R}Hom}}
\newcommand{\usc}[1][m]{\underline{\phantom{#1}}}
\newcommand{\Image}{\operatorname{Image}}
\newcommand{\et}{\op{\acute{e}t}}

\section{Riemann--Hilbert type correspondences and applications to local cohomology invariants}
This section starts with a very rough introduction to some aspects of the
Riemann--Hilbert correspondence (over $\CC$) and, shortly thereafter, to a positive
characteristic analog, developed recently by Emerton and Kisin \cite{EmKis.Fcrys}.
We start at a basic level -- merely motivating the correspondence with the help of a
fundamental example -- but progress quickly to a nontrivial construction which is
central to our applications: the intermediate extension.

The aim is to show how these correspondences can be put to good use in order to
study singularities. Concretely we obtain a new characterization of invariants
arising from local cohomology in terms of \'etale cohomology. One central
interesting aspect of our approach is that the treatment is, up to the use of the
respective correspondence, independent of the characteristic.

The result we will discuss is a description of Lyubeznik's local cohomology
invariants in any characteristic \cite{Lyub.FinChar0} which are originally defined
for a quotient $A=R/I$ of a $n$-dimensional regular local ring $(R,\mathfrak{m})$ to
be
\[
    \lambda_{a,i} \defeq e(H^a_{\mathfrak{m}}(H^{n-i}_I(R)))
\]
where $e(\usc)$ denotes the $D$--module multiplicity (see Section
\ref{manu.loc.inv}). In \cite{Lyub.FinChar0} these are shown to be independent of
the representation of $A$ as the quotient of a regular local ring. In
\cite{BliBon.LocCohomMult} and \cite{Bli.CohomIsol} these invariants were described
in many cases in terms of \'etale cohomology; we discuss these results here as an
application of the aforementioned correspondences:
\begin{thm}\label{thm.main}
Let $A=\Oo_{Y,x}$ for $Y$ a closed $k$--subvariety of a smooth variety $X$. If for
$i \neq d$ the modules $H^{n-i}_{[Y]}(\Oo_X)$ are supported at the point $x$ then
\begin{enumerate}
    \item For $2 \leq a \leq d$ one has
    \[
        \lambda_{a,d}(A)-\delta_{a,d} = \lambda_{0,d-a+1}(A)
    \]
    and all other $\lambda_{a,i}(A)$ vanish.
    \item
    \[
    \lambda_{a,d}(A)-\delta_{a,d} =
    \begin{cases}
    \dim_{\FF_p} H^{d-a+1}_{\{x\}}(Y_{\et},\FF_p) & \text{if $\op{char} k = p$} \\
        \dim_{\CC}   H^{d-a+1}_{\{x\}}(Y_{\op{an}},\CC) & \text{if $k = \CC$}
    \end{cases}
    \]
\end{enumerate}
where $\delta_{a,d}$ is the Kroneker delta function.
\end{thm}
The apparent analogy between the situation over $\CC$ and over $\FF_p$ suggested by
this result is somewhat misleading. As we briefly discuss at the end of this
section, \'etale cohomology with $\FF_p$--coefficients is only a very small part
(the slope zero part) of, say, crystalline cohomology. For crystalline cohomology on
the other hand there are comparison results to singular cohomology (via de Rham
theory) hence Lyubeznik's invariants really capture a different type of information
in characteristic 0 than they do in characteristic $p>0$.

\subsection{Riemann--Hilbert and Emerton--Kisin correspondence}
Let us first fix the following notation: Throughout this section $X$ will be a
smooth scheme over a field $k$ of dimension $n$. Mostly $k$ will be either $\CC$,
the field of complex numbers, or $\FF_q$, the finite field with $q=p^e$ elements, as
in Section 1.

\subsubsection{ The Riemann--Hilbert correspondence}
Now let $k=\CC$. In its simplest incarnation the Riemann--Hilbert correspondence
asserts a one to one map between
\[
    \left\{\begin{matrix}\text{local systems of}  \\ \text{$\CC$-vectorspaces}\end{matrix} \right\} \leftrightarrow \left\{\begin{matrix}\text{locally free coherent $\Oo_X$--modules}\\ \text{with regular singular integrable connection}\end{matrix} \right\}
\]
This correspondence grew out of Hilbert's 21st problem, motivated by work of
Riemann, to find, given a monodromy action at some points, a Fuchsian (regular
singular) differential equation with the prescribed monodromy action at the singular
points.

The correspondence is given via de~Rham theory: A connection is a $k$--linear map
$\nabla: \Mm \to \Omega^1_X \tensor_{\Oo_X} \Mm$ that satisfies the Leibniz rule
$\nabla(rm)=r\nabla(m) + dr \tensor m$. One can extend $\nabla$ in an (up to signs)
obvious way to get a sequence of maps (each denoted also by $\nabla$)
\[
    \Mm \to[\nabla] \Omega^1_X \tensor \Mm \to[\nabla] \Omega_X^2 \tensor \Mm
    \to[\nabla] \ldots \to \Omega^n_X \tensor \Mm
\]
and the connection is called integral if this sequence is a complex (i.e.\
$\nabla^2=0$), called the de~Rham complex $\deR(\Mm)$ associated to the connection
$\nabla$. Now, given such integrable connection its \emph{horizontal sections}
\[
    \Mm^\nabla = \ker \nabla = H^0(\deR(\Mm))
\]
is a local system.

The most trivial, but at the same time most important example, is $\Mm = \Oo_X$ with
$\nabla = d$ the universal differential $\Oo_X \to[d] \Omega^1_X$. This is the
$\Oo_X$--module associated to the system of differential equations
$(\frac{\partial}{\partial x_1},\ldots,\frac{\partial}{\partial x_n}) \cdot f = 0$
where $x_1,\ldots,x_n$ are local coordinates at some point of $X$. The corresponding
local system (the solutions to the differential equation) is of course the constant
local system $\CC$. In fact more is true: By the Poincar\'e Lemma, the (analytic)
de~Rham complex is a resolution of the constant sheaf $\CC$ and hence we can
rephrase this by saying that $\CC$ is quasi-isomorphic to the de~Rham complex.

In view of Grothendieck's philosophy the above correspondence is flawed since the
categories on either side are not closed under any reasonable functors. For example,
the pushforward of a local system is generally not a local system as the inclusion
of a point $\{x\} \to[i] \AA^1$ readily illustrates ($i_*\CC$ is a skyscraper
sheaf). However it is a constructible sheaf (that is one that is locally constant on
each piece of a suitable stratification of $X$), and in fact on constructible
sheaves all the functors one would like to have are defined.

On the other side of the correspondence, modules with integrable connection are
replaced by modules over the ring of differential operators $D_X$. The conditions
one has to impose are holonomicity (which is the crucial finiteness condition) and a
further condition (namely that $\Mm$ is regular singualar) which will not be
considered here. Without giving the precise definition, a holonomic $D_X$--module is
one of minimal possible dimension, and the category they form enjoys a strong
finiteness condition:
\begin{prop}[\cite{Borel.Dmod}]
    The category of holonomic $D_X$--modules is abelian and closed under extensions.
    Every holonomic $D_X$ module has finite length.
\end{prop}
The correspondence should again be given via de~Rham theory. However, it quickly
becomes clear that one has to pass to the derived category. We give the following
simple example as a small indication that the passage to the derived category cannot
be avoided.
\begin{ex}
    Let $X=\AA^1 = \Spec \CC[x]$ such that $D_X=k[x,\frac{\partial}{\partial x}]$.
    Let $\Mm=H^1_{\{0\}}(\Oo_{\AA^1})$ denote the cokernel of the following injection of
    $D_X$--modules
    \[
        \CC[x] \to \CC[x,x^{-1}].
    \]
    Then, as a $\CC$-vectorspace $\Mm$ has a basis consisting of $\frac{1}{x^i}$ for
    $i > 0$. The de~Rham complex associated to $\Mm$ is given
    \[
        \Mm \to[\nabla] \Omega^1_{\AA^1}  \tensor \Mm
    \]
    where the map sends a basis element $\frac{1}{x^{i}} \to dx \tensor
    \frac{-i}{x^{i+1}}$. Hence one immediately sees that $\nabla$ is injective and
    that its cokernel is generated by $dx \tensor \frac{1}{x}$, which is the
    skyscraper sheaf supported at $0$. Hence we have an quasi-isomorphism
    \begin{equation}\label{eq.simpex}
        i_*\CC_{\{0\}} \simeq \deR(H^1_{\{0\}}(\Oo_{\AA^1}))[1].
    \end{equation}
    Note the appearance of the shift $[1]$ which is an indication that one cannot
    avoid to pass to the derived category.
\end{ex}
The ultimate generalization of the basic version above is the Riemann--Hilbert
correspondence as proved by Mebkhout \cite{Mebk.Phd}, Kashiwara
\cite{Kashiwara.FaisConst}, Beilinson and Bernstein \cite{Borel.Dmod}:
\begin{thm}
Let $X$ be a smooth $\CC$--variety. On the level of bounded derived categories,
there is an equivalence between
\[
\left\{\begin{matrix}\text{constructible sheaves of }  \\
\text{$\CC$-vectorspaces}\end{matrix} \right\} \leftrightarrow
\left\{\begin{matrix}\text{holonomic $D_X$--modules}\\ \text{which are regular
singular}\end{matrix} \right\}
\]
The correspondence is given by sending a complex $\Mm^\bullet$ of $D_X$--modules to
$\deR(\Mm^\bullet)= \RHom(\Oo_X,\Mm)$.

This equivalence preserves the \emph{six standard functors} $f^*$ and $f_*$, their
duals $f_!$ and $f^!$ and also $\tensor$ and $\Hom$.
\end{thm}
Via duality on $X$ the de~Rham functor is related to the functor $\Sol(\usc) =
\RHom(\usc,\Oo_X)$. This functor yields therefore an anti-equivalence, and it is
this anti-equivalence which can be obtained in positive characteristic.

\subsubsection{Emerton--Kisin correspondence}
In positive characteristic a naive approach via de~Rham theory does not work due to
the failure of the Poincar\'e lemma. At least one should expect that the
correspondence respects the relationship ``$\FF_p$ corresponds to $\Oo_X$'' in some
way or another. But due to the fact that in characteristic $p$ there are many
functions with derivative zero (namely all $p$th powers), the de~Rham theory is ill
behaved. For the same reason, $D$--module theory in positive characteristic is also
quite different from the one in characteristic zero.

So what is the correct counterpart for constructible sheaves of
$\FF_p$--vector\-spaces? The solution arises from the Artin--Schreyer sequence:
\[
    0 \to \FF_p \to \Oo_X \to[x \mapsto x^p] \Oo_X \to 0.
\]
This is an exact sequence in the \'etale topology. Hence if one views the
Riemann-Hilbert correspondence as a vast generalization of de~Rham theory (in terms
of the  Poincar\'e lemma), then the correspondence of Emerton--Kisin -- to be
outlined shortly -- is an analogous generalization of Artin--Schreyer theory.

Hence, the basic objects one studies in this correspondence is not $D_X$--modules
but rather quasi-coherent $\Oo_X$--modules $\Mm$ which are equipped with an action
of the Frobenius $F$. That is, we have an $\Oo_X$--linear map $F_{\Mm}: \Mm \to
F_*\Mm$. By adjunction, such a map is the same as a map $\theta_\Mm : F^*\Mm \to
\Mm$. Such objects have been studied in various forms for a long time, see for
example \cite{HaSp} or \cite[Exp. XXII]{SGA7.2}.
\begin{defn} An \emph{$\Oo_{X,F}$--module}
$(\Mm,\theta)$ is a quasi-coherent $\Oo_X$--module $\Mm$ together
    with a $\Oo_X$--linear map
    \[
        \theta: F^* \Mm \to \Mm
    \]
    If $\theta$ is an isomorphism $(\Mm,\theta)$ is called \emph{unit}. If is called
    \emph{finitely generated} if $\Mm$ is finitely generated when viewed as a module
    over the non-commutative ring $\Oo_{X,F}$.
\end{defn}
The following two are the essential examples of finitely generated unit
$\Oo_{X,F}$--modules.
\begin{ex}
    For simplicity assume that $X=\Spec R$ is affine. The natural identification
    $F^*\Oo_X \cong \Oo_X$ gives $\Oo_X$ the structure of a finitely generated unit
    $\Oo_{X,F}$--module.

    Let $R_f$ be the localization of $R$ at a single element $f\in R$. The natural
    map
    \[
        F^*R_f = R \tensor_F R_f \to R_f
    \]
    sending $r \tensor \frac{a}{b}$ to $r\frac{a^p}{b^p}$ has a natural inverse
    given by sending $\frac{a}{b} \to ab^{p-1} \tensor \frac{1}{b}$. Hence $R_f$
    is naturally a unit module. In fact, $R_f$ is even finitely generated as a unit
    $R[F]$--module, generated by the element $\frac{1}{f}$, since
    $F(\frac{1}{f})=\frac{1}{f^p}$.
\end{ex}
The main result of \cite{Lyub} about finitely genrated unit $\Oo_{X,F}$--modules,
which makes them a suitable analog of holonomic $D$--modules, is
\begin{thm}
Let $X$ be smooth. In the abelian category of (locally) finitely generated
$\Oo_{X,F}$--modules, every object has finite length.
\end{thm}
\begin{ex}
    Considering a Cech resolution to compute coherent cohomology with support in
some subscheme $Z$, the easy \emph{abelian category} part of the preceding theorem
and the preceding examples imply that $H^i_Z(\Oo_X)$ is naturally a f.g. (finitely
generated) unit $\Oo_{X,F}$--module.
\end{ex}
Now the correspondence that is proven in \cite{EmKis.Fcrys} can be summarized as
follows:
\begin{thm}
For $X$ smooth,  there is an anti-equivalence on the level of derived categories
\[
    \left\{\begin{matrix}\text{constructible sheaves of }  \\ \text{$\FF_p$--vector-spaces on $X_{\et}$}\end{matrix} \right\} \leftrightarrow \left\{\begin{matrix}\text{locally finitely generated}\\ \text{$\Oo_{X,F}$--modules}\end{matrix} \right\}
\]
The correspondence is given by sending a (complex of) $\Oo_{X,F}$--modules
$\Mm^\bullet$ to the constructible sheaf $\Sol(\Mm^\bullet)=\RHom(\Mm,\Oo_X)$ and is
roughly dual to the naive approach of taking the fixed points of the Frobenius
(Artin-Schreyer sequence) alluded to above.

The correspondence preserves certain functors, namely $f^!$, $f_*$ (and $\tensor$)
on the right hand side correspond to $f^*$, $f_!$ on the left hand side.
\end{thm}

\subsubsection{Intermediate extensions}
{F}rom now on we treat the two situations -- characteristic zero regular singular
holonomic $D_X$--modules with the Riemann--Hilbert correspondence on one hand and
positive characteristic finitely generated unit $\Oo_{F,X}$--modules with the
Emerton--Kisin correspondence on the other -- formally as one. The crucial property
they both share is the fact that all modules have finite length. This is key to the
following construction.

Due to the lack of duality in positive characteristic there is no analog of the
functor $j_{!}$ available. However there is an adequate substitute, which still
exists in our context \cite{EmKis.FcrysIntro,Manuel.int}. The substitute we have in
mind is the intermediate extension $j_{!*}$, usually constructed as the image of the
``forget supports map'' from $j_! \to j_*$. This usual definition does not work
since $j_!$ is not available. Nevertheless it turns out that all one needs to define
$j_{!*}$ is the fact that the modules have finite length:
\begin{prop}[\cite{Manuel.int},\cite{EmKis.FcrysIntro}]
    Let $j: U \subseteq X$ be a locally closed immersion of smooth $\FF_p$--schemes (resp.
    $\CC$--schemes). Let $\Mm$ be a finitely generated unit $\Oo_{U,F}$--module.
    Then there is a \emph{unique} submodule $\Nn$ of $R^0j_*\Mm$ minimal with respect
    to the property that $f^!\Nn = \Mm$. This submodule $\Nn$ is called the
    \emph{intermediate extension} and is denoted by $j_{!*}\Mm$.
\end{prop}
\begin{proof}
    The key point is the fact that $\Mm$ has finite length which ensures the existence of
    \emph{minimal} modules with the desired property (any decreasing chain is eventually constant).
    Let $\Nn_1$ and $\Nn_2$ be two modules with the desired property. Since $f^!\Nn_i = \Mm$ their
    intersection cannot be zero. On the other hand the exact sequence
    \[
        0 \to \Nn_1 \cap \Nn_2 \to R^0j_*\Mm \to R^0j_*\Mm/\Nn_1 \oplus R^0j_*\Mm/\Nn_2
    \]
    shows that their intersection is also a finitely generated unit $\Oo_{X,F}$--module (resp. r.s. holonomic $D_X$--module)
    as it is the kernel of a map of such modules. By minimality one has $\Nn_1 = \Nn_2$
    showing uniqueness.
\end{proof}
Under the correspondence the intermediate extensions behave well: In the situation
as above we have
\[
    \Sol(j_{!*}\Mm) = \Image(j_!\Sol(\Mm) \to R^0j_*\Sol(\Mm))
\]
such that they do in fact correspond to the intermediate extensions on the
constructible side, where they can be identified as the image of the ``forget
supports'' map.

For the basic computations that follow we list some key properties both
correspondences enjoy:
\begin{enumerate}
    \item \[
                \Sol(\Oo_X) = \begin{cases} \FF_p[n] &\text{if $k=\FF_p$} \\ \CC[n] &\text{if
                $k=\CC$}\end{cases}
            \]
    \item There are functors $f^!$ and $f_*$ which behave under the correspondence in the expected way.
    \item For $Y \subseteq X$ a subvariety, local cohomology is defined in the categories and satisfies the
    triangle (a highbrow way of writing the long exact sequence for cohomology with
    supports)
        \[
            \R\Gamma_{[Y]}\Mm^\bullet \to \Mm^\bullet \to \R j_* j^*
            \Mm^\bullet \to[+1]
        \]
        where $j$ is the open inclusion of the complement of $Y$ into $X$.
    \item Again, let $i: Y \hookrightarrow X$ be the inclusion of a closed subset,
    then
    \[
        \Sol \circ \R\Gamma_{[Y]} \cong \R i_!i^* \circ \Sol.
    \]
    This follows via the preceding two items and the triangle $\R
    j_!j^!\Ll \to \Ll \to \R i_!i^{-1} \Ll \to[+1]\quad$ for a (complex of) constructible sheaves
    $\Ll$.
    \item The preceding items allow us to compute (see also equation
    (\ref{eq.simpex}) on page \pageref{eq.simpex})
    \[
        \Sol(\R\Gamma_Y(\Oo_X)) = i_!i^* \FF_p[n] = i_!\FF_p|_Y[n]
    \]
    in positive characteristic and respectively $i_!\CC|_Y[n]$ in characteristic zero.
\end{enumerate}

\subsection{Lyubeznik's local cohomology invariants}\label{manu.loc.inv}
Let $A=R/I$ for $I$ an ideal in a regular (local) ring $(R,m)$ of dimension $n$ and
containing a field $k$. The main results of \cite{Lyub.FinChar0,HuSha.LocCohom}
state that the local cohomology module $H^a_m(H^{n-i}_I(R))$ is injective and
supported at $m$. Therefore it is a finite direct sum of $e=e(H^a_m(H^{n-i}_I(R)))$
many copies of the injective hull $E_{R/m}$ of the residue field of $R$. Lyubeznik
shows in \cite{Lyub.FinChar0} that this number
        \[
        \lambda_{a,i}(A) \defeq e(H^a_m(H^{n-i}_I(R)))
        \]
does not depend on the auxiliary choice of $R$ and $I$. At the same time, this
number $e(\Mm)$ is the multiplicity of the holonomic $D_X$--module $\Mm$,
respectively the finitely generated unit $\Oo_{X,F}$--module $\Mm$.

If $A$ is a complete intersection, these invariants are essentially trivial (all are
zero except $\lambda_{d,d}=1$ where $d=\dim A$). In general $\lambda_{a,i}$ can only
be nonzero in the range $0 \leq a,i \leq d$. These invariants were first introduced
by Lyubeznik in \cite{Lyub.FinChar0} and further studied by Walther in \cite{Walt}.
In \cite{LopezSabbah} Garcia-L\'opez and Sabbah show Theorem \ref{thm.main} in the
case of an isolated complex singularity.

We now indicate briefly the proof of \ref{thm.main}. As shown in
\cite{BliBon.LocCohomMult} the condition imposed on the singularities in the Theorem
\ref{thm.main} easily implies (via the spectral sequence $E_2^{a,j} =
H^a_{[x]}H^{j}_{[Y]}(\Oo_X) \Rightarrow H^{a+j}_{[x]}(\Oo_X)$) part one. Part two is
the point where the correspondences enter into the picture: The idea is of course to
use that
\[
    e(H^a_m(H^{n-d}_{Y}(\Oo_X))) = \dim \Sol(H^a_{\{x\}}(H^{n-d}_Y(\Oo_X)))
\]
by the correspondence, and then to compute the right hand side. As it is written
here the right hand side is however not computable. The trick now is to replace
$H^{n-d}_I(R)$ by $j_{!*} H^{n-d}_Y(\Oo_X)|_{X-\{x\}}$, which is easily checked to
not affect our computation (long exact sequence for $\Gamma_{\{x\}}$). Now, the
assumption on the singularity that for $i\neq d$ the module $H^{n-i}_Y(\Oo_X)$ is
supported at $x$ can simply be rephrased as $H^{n-d}_Y(\Oo_X)|_{X-\{x\}}=
\R\Gamma_{Y-\{x\}}(\Oo_{X-\{x\}})[n-d]$.

Using that $\Sol$ commutes with $j_{!*}$ and the fact that
\[
    \Sol(\R\Gamma_{Y-\{x\}}(\Oo_{X-\{x\}})) = i_!(\FF_p)_{Y-\{x\}}[n]
\]
one obtains that
\[
    e(H^a_m(H^{n-d}_{Y}(\Oo_X))) = \dim (H^{-a}j_{!*}i_!(\FF_p)_{Y-\{x\}}[d]).
\]
To compute the right hand side is now a feasible task  that yields the desired
result (feasible due to the fact that $j$ is just the inclusion of the complement of
a point, which makes it possible to effectively understand and calculate $j_{!*}$;
see \cite{BliBon.LocCohomMult,Bli.CohomIsol} for details).

\subsection{Comparison via crystalline cohomology} We close this
section with some remarks regarding the behaviour of these invariants under
reduction to positive characteristic. There are by now classical examples that show
that local cohomology does not behave well under reduction so one would expect that
the invariants $\lambda_{a,i}$ do not behave well either. On a superficial level,
glancing at Theorem \ref{thm.main}, one might however suspect a complete analogy
between positive and zero characteristic.

However, this is not true. The difference stems from the difference between the
cohomology theories which describe $\lambda_{a,i}$. In positive characteristic, this
is \'etale cohomology with coefficients in $\FF_p$, in characteristic zero however
it is (topological) cohomology. Under reduction mod $p$ the former only constitutes
a very small part of the latter, namely the part of slope zero.

For simplicity consider the situation of $Y \subseteq \PP^n$ a smooth projective
variety. Now. the local ring $A$ of the cone of this projective embedding has an
isolated singularity at its vertex and one can study the invariants
$\lambda_{a,i}(A)$. Theorem \ref{thm.main} shows that $\lambda_{a,d}(A)$ are
described (excision) by $H^{d-a}(Y_{\et},\FF_p)$ and $H^{d-a}(Y_{\op{an}},\CC)$ in
positive and zero characteristic respectively. So far everything appears in complete
analogy.

But let us now consider reduction mod $p$ and let $Y$ be the special fiber of the
smooth family $\Yy \to \Spec W(k)$, where $W(k)$ is the ring of Witt vectors over
$k$ and $K$ denotes its field of fractions. Via the comparison results between
topological cohomology, de~Rham cohomology and Berthelot's crystalline cohomology
(see \cite{BertOgus})
\[
    H^{d-a}(\Yy_\CC,\CC) \to[\cong] H^{d-a}_{\deR}(\Yy_K) \cong H^{d-a}(Y/K)
\]
it turns out that $H^{d-a}(Y_{\et},\FF_p)$ is only a very small part of the
crystalline cohomology $H^{d-a}(Y/K)$, namely the part on which the Frobenius acts
with eigenvalue zero. This is an even smaller part then the part $H^{d-a}(Y/K)^{<1}$
which was of great importance in section \ref{subsec.16}. From this point of view it
is not surprising to find the characteristic zero side to capture much more
information than the positive characteristic side. Even though the positive
characteristic side is therefore lacking the topological insight, the information
one obtains is still very valuable. For example it is used very effectively to study
$L$--functions in \cite{EmKis.Fcrys} and also \cite{BoPi}, which uses a similar
correspondence for that purpose.




\begin{thebibliography}{[KLR73]}
\bibitem{Ax} Ax, J.: Zeroes of polynomials over finite fields,
Amer. J. Math. {\bf 86} (1964), 255-261.
\bibitem{BBE} Berthelot, P.,  Bloch, S., Esnault, H.: On Witt vector cohomology
for singular varieties, preprint 2005, 42 pages.
\bibitem{BertOgus}
Berthelot, P., Ogus, A.: Notes on crystalline cohomology. Princeton University
Press, Princeton, N.J.; University of Tokyo Press, Tokyo, 1978. vi+243 pp.
\bibitem{BliBon.LocCohomMult}
Blickle, M.,   Bondu, R.: Local cohomology multiplicities in
  terms of local cohomology, 2004, to appear in Ann. Inst. Fourier.
\bibitem{Bli.CohomIsol}
Blickle, M. : Lyubeznik's invariants for cohomologically isolated singularities,
preprint 2005.
\bibitem{Manuel.int}
Blickle, M.: The intersection homology {$D$}--module in finite
  characteristic, Math. Ann. \textbf{328} (2004), 425--450.
\bibitem{Bl} Bloch, S.: Algebraic $K$-theory and crystalline
cohomology, Publ. Math. I.H.\'E.S. {\bf 47} (1977), 187--268.
\bibitem{B} Bloch, S.: Lectures on Algebraic cycles, Duke
University Mathematics Series, IV, (1980).
\bibitem{Bl86} Bloch, S.: Algebraic cycles and higher K-theory, Adv. in Math. {\bf 61}, No.3 (1986), 267-304.
\bibitem{BlEs03a} Bloch, S., Esnault, H.: An additive version of higher Chow groups,  Ann. Sci. \'Ecole Norm. Sup. 4
 {\bf 36},No. 3 (2003), 463-477.
\bibitem{BlEs03b} Bloch, S., Esnault, H.: The additive dilogarithm, Doc. Math., J. DMV Extra Vol. (2003), 131-155.
\bibitem{BEL}  Bloch, S., Esnault, H.,   Levine, M.: Decomposition of the
diagonal and eigenvalues of Frobenius for Fano hypersurfaces,
 Am. J. of  Mathematics,  {\bf 127} no.1 (2005), 193-207.
\bibitem{BoPi}
Böckle, G., Pink, R.: Cohomological {T}heory of crystals
  over function fields, uncirculated preprint, 2005.
\bibitem{Borel.Icohom}
Borel, A. et al: Intersection cohomology, Progress in Mathematics,
  vol.~50, Birkh\"auser Boston Inc., Boston, MA, 1984, Notes on the seminar
  held at the University of Bern, Bern, 1983, Swiss Seminars.
\bibitem{Borel.Dmod}
Borel, A.,  Grivel, P.-P,  Kaup, B., Haefliger, A., Malgrange, B., and Ehlers, E.:
  Algebraic ${D}$-modules, Academic Press Inc., Boston, MA, 1987.
\bibitem{dJ}  de Jong, A. J.: Smoothness, semi-stability and alterations, Publ. Math. I.H.\'E.S. {\bf 83} (1996), 51-93.
\bibitem{DHodge} Deligne, P.: Th\'eorie de Hodge II, Th\'eorie de Hodge III,  Publ. Math. I.H.\'E.S.
{\bf 40} (1971), 5--57, {\bf 44} (1974), 5--77.
\bibitem{DWeil} Deligne, P.: La conjecture de Weil. I, La conjecture de Weil. II,
Publ. Math. I.H.\'E.S.  {\bf 43} (1974), 273--307, {\bf 52}, (1980) 137--252.
\bibitem{Dint} Deligne, P.: Th\'eor\`eme d'int\'egralit\'e, Appendix  to
 Katz, N.: Le niveau de la cohomologie des intersections compl\`etes, Expos\'e XXI
in SGA 7, Lect. Notes Math. vol. {\bf 340}, 363-400, Berlin Heidelberg New York
Springer 1973.
\bibitem{DE} Deligne, P., Esnault, H.: Appendix to
``Deligne's integrality theorem in unequal characteristic and rational points over
finite fields'', preprint 2004, 5 pages, appears in the Annals of Mathematics.
\bibitem{SGA7.2} Deligne, P., Katz, N.: Seminar de Geometrie Algebrique, SGA 7.II, Lect. Notes Math. vol. {\bf 340}
Spinger 1973, Berlin Heidelberg New York.
\bibitem{EmKis.Fcrys}
Emerton, M.,   Kisin, M.: Riemann--Hilbert correspondence for unit
  $\mathcal{F}$-crystals, {Ast�isque} \textbf{293} (2004), vi+257 pp.
\bibitem{EmKis.FcrysIntro}
Emerton, M.,  Kisin, M.: An introduction to the
  Riemann--Hilbert correspondence for unit $\mathcal{F}$-crystals.,
  Geometric aspects of Dwork theory. Vol. I, II,  677--700,
Walter de Gruyter GmbH \& Co. KG,
  Berlin, 2004.
\bibitem{EHodgetyp} Esnault, H.: Hodge type of subvarieties of $ \P^n$  of small degrees, Math. Ann.
{\bf 288} (1990), 549 - 551.
\bibitem{ENS} Esnault, H., Nori, M.,  Srinivas, V.: Hodge type of projective varieties of
low degree. Math. Ann. {\bf 293} (1992), 1-6.
\bibitem{EFano} Esnault, H. : Varieties over a finite field with trivial Chow group of 0-cycles
have a rational point, Invent. math. {\bf 151} (2003), 187-191.
\bibitem{Eeigv} Esnault, H.: Eigenvalues of Frobenius acting on the $\ell$-adic cohomology
of complete intersections of low degree,  C. R. Acad. Sci. Paris, Ser. I {\bf 337}
(2003), 317-320.
\bibitem{EK}  Esnault, H., Katz. N: Cohomological divisibility and point count
divisibility, Compositio Mathematica {\bf 141} 01 (2005), 93-100.
\bibitem{EappFakh} Esnault, H.: Appendix to `` Congruences for rational points on varieties over
finite  fields'' by N. Fakhruddin and C. S. Rajan,
Math. Ann. {\bf 333} (2005), 811-814.
\bibitem{Eunqchar} Esnault, H.: Deligne's integrality theorem in unequal characteristic and rational
points over finite fields, preprint 2004, 10 pages, appears in the Annals of
Mathematics.
\bibitem{FR} Fakhruddin, N., Rajan, C. S.: Congruences for rational points on varieties over finite  fields,
Math. Ann. {\bf 333} (2005), 797-809.
\bibitem{Fu} Fujiwara, K.: A Proof of the Absolute Purity Conjecture (after Gabber), in Algebraic Geometry 2000, Azumino, Advanced Studies in Pure Mathematics {\bf 36} (2002), Mathematical Society of Japan, 153-183.
\bibitem{LopezSabbah}
Garcia-L\`opez, R.,  Sabbah, C.: Topological computation of
  local cohomology multiplicities, Collect. Math. \textbf{49} (1998), no.~2-3,
  317--324, Dedicated to the memory of Fernando Serrano.
\bibitem{Gr} Grothendieck, A.: Formule de Lefschetz et rationalit\'e
des fonctions $L$, S\'eminaire Bourbaki {\bf 279}, 17-i\`eme ann\'ee (1964/1965),
1-15.
\bibitem{HaSp}
Hartshorne, R., Speiser, R.: Local cohomological dimension in
  characteristic $p$, Annals of Mathematics \textbf{105} (1977), 45--79.
\bibitem{HeMa01} Hesselholt, L., Madsen, I.: On the $K$-theory of nilpotent endomorphisms,
Homotopy methods in algebraic topology (Boulder, CO, 1999), Contemp. Math. {\bf
271}, Amer. Math. Soc., Providence, RI (2001), 127-140.
\bibitem{HuSha.LocCohom}
Huneke, C. L., Sharp, Ro.Y.: Bass numbers of local cohomology
  modules, Trans. Amer. Math. Soc. \textbf{339} (1993), no.~2, 765--779.
\bibitem{I} Illusie, L.: Complexe de de Rham-Witt et cohomologie
cristalline, Ann. \'Ec. Norm. Sup. 4  {\bf 12} (1979), 501--661.
\bibitem{Kashiwara.FaisConst}
Kashiwara, M.: Faisceaux constructibles et syst\`emes holon\^omes
  d'\'equations aux d\'eriv\'ees partielles lin\'eaires \`a points singuliers
  r\'eguliers, S\'eminaire Goulaouic-Schwartz, 1979--1980, \'Ecole
  Polytech., Palaiseau, 1980, pp.~Exp. No. 19, 7.
\bibitem{Ka}  Katz, N.:  On a theorem of Ax,  Amer. J. Math.
{\bf 93}  (1971),  485-499.
\bibitem{Ko} Koll\'ar, J.:
Rational curves on algebraic varieties,
 Ergebnisse der Mathematik
und ihrer Grenzgebiete. 3. Folge, {\bf 32} (1996), Springer-Verlag, Berlin, 1996.
\bibitem{L} Lang, S.: Cyclotomic points, very anti-canonical
varieties, and quasi-algebraic closure, preprint 2000.
\bibitem{Lyub.FinChar0}
Lyubeznik, G.: Finiteness properties of local cohomology modules (an
  application of ${D}$-modules to commutative algebra), Invent. math.
  \textbf{113} (1993), no.~1, 41--55.
\bibitem{Lyub}
Lyubeznik, G.: $\mathcal{F}$-modules: an application to local
  cohomology and {$D$}-modules in characteristic $p>0$, Journal f{\"ur} reine
  und angewandte Mathematik \textbf{491} (1997), 65--130.
\bibitem{Ma} Manin, Yu.: Notes on the arithmetic of Fano
threefolds, Compos. math. {\bf 85} (1993), 37-55.
\bibitem{Mebk.Phd}
Mebkhout, Z.: Cohomologie locale des espaces analytiques complexes, Ph.D.
  thesis, Universit{\'e} Paris {V}{I}{I}, 1979.
\bibitem{NeSu89} Nesterenko, Y. P., Suslin, A.: Homology of the general linear group over a local ring, and Milnor's K-theory, Math. USSR-Izv. {\bf 34}, No.1 (1990), 121-145.
\bibitem{R} R\"ulling, K.: The generalized de Rham-Witt complex over a field is a complex of zero-cycles, preprint 2005.
\bibitem{Se} Serre, J.-P.: Sur la topologie des vari\'et\'es
alg\'ebriques en caract\'eristique $p$, in \textit{Symposium internacional de
topolog\'{\i}a algebraica} (International symposium on algebraic topology), 24--53,
Universidad Nacional Autonoma de Mexico and UNESCO, Mexico City, 1958.
\bibitem{To92} Totaro, B.: Milnor K-theory is the simplest part of algebraic K-theory, K-theory {\bf 6}, No.2 (1992), 177-189.
\bibitem{V02} Voevodsky, V.: Motivic cohomology groups are isomorphic to higher Chow groups in any characteristic,
Int. Math. Res. Not. {\bf 7} (2002), 7, 351-355.
\bibitem{Walt} Walther, U.: On the {L}yubeznik numbers of a local ring, Proc. Amer. Math. Soc.
\textbf{129} (2001), No. 6, 1631--1634 (electronic).
\bibitem{Wi36} Witt, E.: Zyklische K\"orper und Algebren der Charakteristik $p$ vom Grad $p^n$. Struktur diskret bewerteter perfekter K\"orper mit vollkommenem Restklassenk\"orper der Charakteristik $p$,
 J. Reine Angew. Math. {\bf 176} (1936), 126-140.







\end{thebibliography}
\end{document}